\documentclass[12pt]{article}
\usepackage{epsfig}
\usepackage{amssymb}

\newtheorem{lem}{Lemma}
\newtheorem{thm}{Theorem}

\begin{document}
\centerline{\large\bf Analytic proof of the existence of the Lorenz attractor}
\centerline{\large\bf in the extended Lorenz model}
\medskip
\centerline{\bf Ovsyannikov I.I.$^{1,2}$, Turaev D.V.$^{1,3}$}
\medskip
\centerline{\it $^1$ Lobachevsky State University of Nizhny Novgorod,}
\centerline{\it $^2$ Universit\"at Bremen,}
\centerline{\it  Jacobs University Bremen,}
\centerline{\it $^3$ Imperial College London,}
\centerline{\it Joseph Meyerhoff Visiting Professor, Weizmann Institute of Science}
\centerline{\it E-mail: ivan.i.ovsyannikov@gmail.com, d.turaev@imperial.ac.uk}
\bigskip

\abstract{We give an analytic (free of computer assistance) proof of the existence of a classical
Lorenz attractor for an open set of parameter values of the Lorenz model in the form of 
Yudovich-Morioka-Shimizu. The proof is based on detection of a homoclinic butterfly
with a zero saddle value and rigorous verification of one of the Shilnikov criteria for the birth
of the Lorenz attractor; we also supply a proof for this criterion. The results are applied in order to give 
an analytic proof of the existence of a robust, pseudohyperbolic strange attractor (the so-called discrete 
Lorenz attractor) for an open set of parameter values in a $4$-parameter family of three-dimensional Henon-like diffeomorphisms.}



\section{Introduction}

The main goal of the paper is a proof of the birth of the Lorenz attractor in the following model:
\begin{equation}
\begin{array}{l}
\dot X = Y, \\
\dot Y = - \lambda Y + \gamma X (1 - Z) - \delta X^3 \\
\dot Z = - \alpha Z + \beta X^2,
\end{array}
\label{ELorenz}
\end{equation}
where parameters $\alpha$, $\beta$, $\gamma$, $\delta$ and $\lambda$ can take arbitrary finite values. It is well-known \cite{Y78, SM78, bel84} that the classical Lorenz equations
\begin{equation}
\begin{array}{l}
\dot x = -\sigma (x - y), \\
\dot y = r x - y - x z \\
\dot z = - b z + x y
\end{array}
\label{Lorenz}
\end{equation}
can be always brought to form (\ref{ELorenz}) by the transformation $\displaystyle x = X \sqrt{2}, \; y = \sqrt{2} \left(\frac{Y}{\sigma} + X \right), \;
z = (r - 1) Z + \frac{X^2}{\sigma}$ with $\displaystyle \alpha = b, \;\; \beta = \frac{2 \sigma - b}{\sigma (r - 1)}, \;\; \gamma = \sigma (r - 1), \;\;
\delta = 1, \;\; \lambda = \sigma + 1$. We therefore call system (\ref{ELorenz}) the {\em extended Lorenz model}.

Not all values of parameters in system (\ref{ELorenz}) correspond to real values of the parameters $(r, \sigma, b)$, so the proof of the Lorenz attractor
in model (\ref{ELorenz}) may not always imply the existence of the Lorenz attractor in the classical 
Lorenz model. In particular, the values of $\alpha$ and $\lambda$,
for which the existence of the Lorenz attractor is proved in the present paper (see Theorem~\ref{thm:main}), do not correspond to the Lorenz model;
instead, for these parameters values one can transform system (\ref{ELorenz}) to the system obtained from the Lorenz model by the time reversal.

Nevertheless, system (\ref{ELorenz})
is interesting by itself because it is a normal form for certain codimension-3 bifurcations of equilibria and periodic
orbits \cite{SST93,PST98}. Thus, our result on the existence of the Lorenz attractor in system (\ref{ELorenz}) automatically proves the emergence of
the Lorenz attractor (or its discrete analogue) in a class of such bifurcations and, hence, in a large set of systems of various nature.
In particular, in this paper, using Theorem~\ref{thm:main} we prove the existence of discrete Lorenz attractors in a class of three-dimensional
polynomial maps (3D Henon maps).

We consider the case $\gamma > 0$, so by scaling the time and the coordinates we can make $\gamma = 1$.
When $\alpha \beta > 0$, with the scaling $\beta = \alpha \delta / B$ one can bring system (\ref{ELorenz}) to the form:
\begin{equation}
\begin{array}{l}
\dot X = Y, \\
\dot Y = X - \lambda Y - X Z - B X^3 \\
\dot Z = - \alpha (Z - X^2).
\end{array}
\label{nfSST}
\end{equation}
In \cite{SST93} it was shown numerically that system (\ref{nfSST}) possesses
the Lorenz attractor in an open set of values of $(\alpha, \lambda)$ for each $B \in (-1/3, +\infty)$.
In this paper we obtain an analytic proof of the same fact
for sufficiently large $B$. The idea is that for $B = +\infty$ equations (\ref{nfSST}) can be solved explicitly\footnote{The same holds for $B = -1/3$, which case will be considered in a forthcoming paper; the case of arbitrary $B > -1/3$ remains out of reach.}.
Therefore for large $B$ the system can be analysed using asymptotic expansions.

We perform such expansion for system (\ref{ELorenz}) with $\gamma >0$,  $\delta >0$. By coordinates and time scaling it can be brought
to the form:
\begin{equation}
\begin{array}{l}
\dot X = Y, \\
\dot Y = X - \lambda Y - X Z - X^3 \\
\dot Z = - \alpha Z + \beta X^2.
\end{array}
\label{nfY}
\end{equation}
Obviously, large $B$ in system (\ref{nfSST}) corresponds to small $\beta$ in (\ref{nfY}).
The following theorem is our main result.

\begin{thm}
\label{thm:main}
For each sufficiently small $\beta$, in the $(\alpha, \lambda)$--plane there exists a point $\alpha = \alpha(\beta), \lambda=\lambda(\beta)$
for which system $($\ref{nfY}$)$
possesses a double homoclinic loop $($a homoclinic butterfly$)$ to a saddle equilibrium with zero saddle value.
This point belongs to the closure of a domain $V_{LA}(\beta)$ in the $(\alpha, \lambda)$--plane for which system
$($\ref{nfY}$)$ possesses an orientable Lorenz attractor.
\end{thm}

In \S~\ref{sect:crit} we give a definition of the Lorenz attractor, following \cite{ABS82}. This definition
requires a fulfilment of a version of a cone condition at each point of a certain absorbing domain. We verify this
condition by using one of the Shilnikov criteria \cite{Sh81, Sh02}. A partial case of these criteria
was proven by Robinson \cite{R89,R92}, however it was done under additional ``smoothness of the foliation'' assumption which is not fulfilled by system (\ref{nfY}) at small $\beta$.
Therefore we provide a full proof of the part of Shilnikov criteria relevant to our situation, see Theorem~\ref{Shil}.
In \S~\ref{sect:main} we prove Theorem~\ref{thm:main}.

In \S~\ref{sect:nf} we prove (Theorem~\ref{thmHen}) the existence of discrete Lorenz attractors in a class of three-dimensional Henon maps.
A discrete Lorenz attractor is a generalization of the strange attractor in the period map of a time-periodic perturbation of an autonomous flow with
a Lorenz attractor \cite{GOST05, GGOT13}. It is a pseudo-hyperbolic attractor in the sense of \cite{ST98, ST08}; like in the continuous-time
Lorenz attractor, each orbit in the discrete Lorenz attractor has positive maximal Lyapunov exponent, and this property is robust
with respect to small smooth perturbations of the map. 

The discrete Lorenz attractors were found numerically in several 3D Henon maps
\cite{GOST05, GGOT13, GGKT14} and in certain models of non-holonomic mechanics \cite{GG15, GGK13}. The 3D Henon maps are particularly
important because they serve as a zeroth order approximation to the rescaled first-return maps near various types of homoclinic and heteroclinic tangencies \cite{GTS93c, GTS08, GTS09, GMO06, GO13, GOT14}. Therefore, because of the robustness of the discrete Lorenz attractor,
by showing that a 3D Henon map has such attractor we also show that the corresponding homoclinic or heteroclinic
tangency bifurcation produces such attractor; moreover, when the corresponding tangency is persistent (e.g. due to
the Newhouse mechanism \cite{N77,PV94,GST97, GST93b}) one obtains infinitely many coexisting discrete Lorenz attractors in a generic map from the corresponding Newhouse domain (see \cite{GO13} for examples and more explanation). Obviously, in order to implement this construction completely
rigorously, the numerical evidence is not enough, i.e. an analytic proof of the existence of the discrete Lorenz attractor in
various classes of Henon maps is needed. Our Theorem~\ref{thmHen} in \S~\ref{sect:nf} is the first example of such proof.

We find that a certain codimension-4 bifurcation happens at certain parameter values to the fixed point of a 3D Henon map and
show that the normal form for this bifurcation is given by the period map of a small time-periodic perturbation of
system (\ref{nfY}) at small $\beta$. Then, Theorem~\ref{thm:main} implies Theorem~\ref{thmHen} immediately.
A less degenerate bifurcation (of codimension 3) that also occurs in the same class of maps produces, as a normal form,
system (\ref{nfSST}) with $B = 0$. This particular case is known as the Morioka-Shimizu system \cite{SM80}.
While a very detailed numerical analysis of the dynamics of this model is available \cite{ASh91,ASh92,Ru92,Ru93,ASh14} and the existence
of the Lorenz attractor in the Morioka-Shimizu system causes no doubt, a rigorous analytic proof of this fact is lacking (see discussion in \cite{TT11}).
Its absence impedes the further development of the mathematical theory of the discrete Lorenz attractors. Obtaining such proof
is a challenging and important problem; in particular, it would be interesting to know if a rigorous numerical proof 
of the existence of the Lorenz attractor in the Morioka-Shimizu system is possible, similar to that reported in \cite{Tucker0,Tucker} for the
classical Lorenz model.

\section{Shilnikov criterion of the birth of Lorenz attractor} \label{sect:crit}

\begin{figure} [htb]
\centerline{\epsfig{file=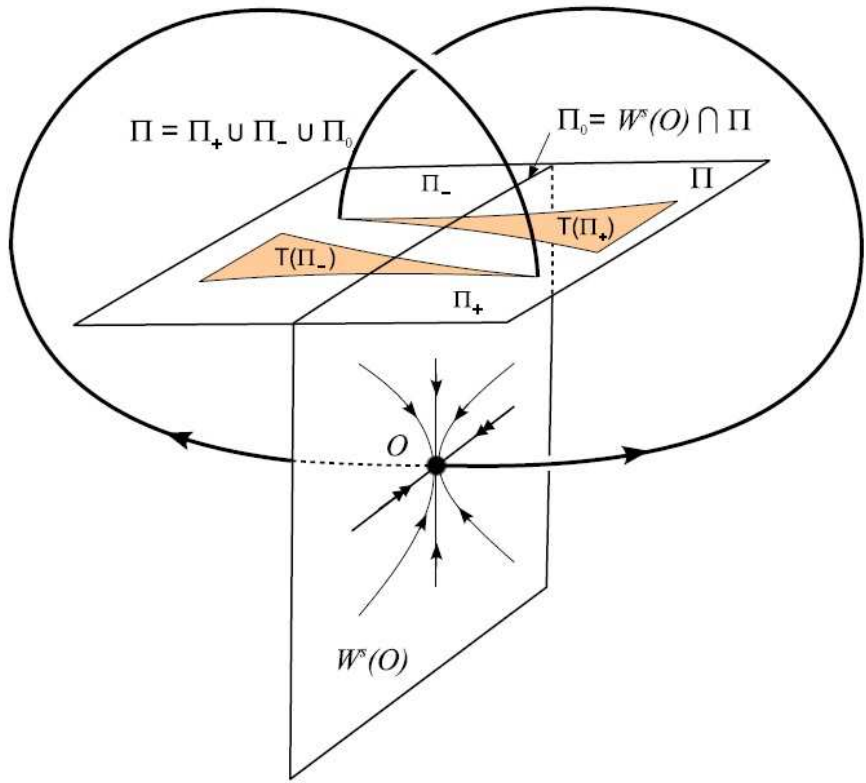, width=12cm}}
\caption{{\footnotesize Afraimovich-Bykov-Shilnikov geometric Lorenz model.}} \label{fig:lorsec}
\end{figure}

Let us recall a definition of the Lorenz attractor. There are several approaches to it, which include classical geometric models by
Guckenheimer-Williams \cite{G76,GW} and Afraimovich-Bykov-Shilnikov \cite{ABS77, ABS82} and their later generalisations
\cite{MPP,AP}. Here we define the Lorenz attractor in the Afraimovich-Bykov-Shilnikov sense.

Consider a smooth system of differential equations which has a saddle equilibrium state $O$
with one-dimensional unstable manifold. Take a cross-section $\Pi$ transverse to a piece of the stable manifold
$W^s(O)$ (for system (\ref{nfY}) the cross-section will be the plane $z=\beta/\alpha$).
The one-dimensional manifold $W^u(O)$ consists of three orbits: $O$ itself, and two separatrices, $\Gamma_+$ and $\Gamma_-$.
Let both $\Gamma_+$ and $\Gamma_-$ intersect the cross-section at some points $M_+$ and $M_-$. Let $\Pi_0$ be
the first intersection of the cross-section with $W^s(O)$, and $\Pi_+$ and $\Pi_-$ be the regions on the opposite sides from $\Pi_0$, so $\Pi=\Pi_+\cup\Pi_-\cup \Pi_0$. Assume that
in the cross-section there is a region $D\subset \Pi_0$ which is forward-invariant with respect to the Poincare map $T$ defined by the orbits of the system.
When a point $M$ approaches $\Pi_0$ from the side of $\Pi_+$ we have $T M \to M_+$,
and when it tends to $\Pi_0$ from the side of $\Pi_-$ we have $T M \to M_-$ (see Fig. \ref{fig:lorsec}). Thus, the invariance of $D$ implies that $M_+\in D$ and $M_-\in D$.
Let us introduce coordinates $(u,v)$
in $\Pi$ such that $\Pi_0$ has equation $u=0$, $\Pi_+$ corresponds to $u>0$, and $\Pi_-$ corresponds to $u<0$. We write the map $T$ in the form
$$\bar u=f_\pm(u,v), \qquad \bar v=g_\pm(u,v),$$
where ``$+$'' corresponds to $u>0$ and ``$-$'' corresponds to $u<0$. The map is smooth
at $u\neq 0$ but it becomes singular as $u\to 0$ because the return time to the
cross-section tends to infinity as the starting point approaches
the stable manifold of the equilibrium $O$. One can analyze the behaviour of the orbits that pass
a small neighbourhood of $O$ by solving Shilnikov boundary value problem (see e.g. formulas from \cite{book}, Section 13.8).
The result gives us the following asymptotics
\begin{equation}\label{pmps}
f_\pm = u_\pm\pm A_\pm |u|^\nu + o(|u|^\nu),\qquad g_\pm = v_\pm\pm B_\pm |u|^\nu + o(|u|^\nu),
\end{equation}
where $(u_+,v_+) = M_+$ and $(u_-,v_-) = M_-$, and $\nu=|\lambda_1|/\gamma$ where $\gamma$ is the positive eigenvalue of the linearisation matrix at $O$ and $\lambda_1$
is the eigenvalue which is nearest to the imaginary axis among the eigenvalues with negative real part; it is assumed that $\lambda_1$ is real and the rest of the eigenvalues
lies strictly farther from the imaginary axis. We assume that the saddle value $\lambda_1+\gamma$ is positive, i.e. $\nu<1$. This implies that if the {\em separatrix values}
$A_\pm$ are non-zero then $\displaystyle \frac{\partial f_\pm}{\partial u}\to\infty$, $\displaystyle \frac{\partial g_\pm}{\partial v}\to 0$
as $u\to 0$, i.e. the map near $\Pi_0$ is expanding in the $u$-direction and contracting in the $v$-directions. The main assumption of the Afraimovich-Bykov-Shilnikov model
\cite{ABS82, ABS77} is that this hyperbolicty property extends from the neighbourhood of $\Pi_0$ to the entire invariant region $D$.
For the first return  map in the form (\ref{pmps})  this assumption is written as the following set of inequalities:
\begin{equation}\label{abscon}
\begin{array}{l} \displaystyle
\|(f'_u)^{-1}\|_{_\circ} < 1,\qquad  \|g'_v\|_{_\circ} < 1,\\ \\ \displaystyle
\|g'_u (f'_u)^{-1}\|_{_\circ} \cdot \|f'_v\|_{_\circ} < (1-\|(f'_u)^{-1}\|_\circ)(1-\|g'_v\|_{_\circ}),
\end{array}
\end{equation}
where we omitted the indices $\pm$; the notation $\|\cdot\|_\circ = \sup_{(u,v)\in D\cap \Pi_\pm} \|\cdot\|$ is used.
These inequalities are an equivalent form of
the so-called invariant cone conditions, which were later used in \cite{BuSi, LuMe} for studying statistical properties of the Lorenz attractor,
in \cite{Vul, Tucker} for the numerical verification of the match between the geometrical Lorenz model and the Lorenz system itself, etc.

In \cite{ABS82} it was shown that conditions (\ref{abscon}) imply the existence of a stable invariant foliation $F^{ss}$ of $D$ by leaves
of the form $u=h(v)$ where $h$ is a certain uniformly Lipshitz function (different leaves are given by different functions $h$).
One can, in fact, show that the foliation $F^{ss}$ is $C^1$-smooth, see \cite{Ro81, ShaShi93}; however, this is not important here.
The map $T$
takes each leaf of the foliation into a leaf of the same foliation and is uniformly contracting on each of them. The distance between the leaves
grows with the iterations of $T$, i.e. the quotient map $\hat T$
obtained by the factorisation of $T$ over the invariant foliation is an expansive map of the interval
obtained by the factorisation of $D$ over the foliation.
This map has a discontinuity at $u = 0$: $\hat T(+0) = u_+$, $\hat T(-0) = u_-$ (see Fig. \ref{fig:lorsec}).
The fact that the quotient map is expansive makes an exhaustive analysis of its dynamics possible, see \cite{Mal89, Mal91}.
The dynamics of the original map $T$ and the corresponding suspension flow can be recovered from that of $\hat T$.

The structure of the Lorenz attractor is described by the following theorem due to Afraimovich-Bykov-Shilnikov.\\

{\bf Theorem} \cite{ABS82, ABS77}
{\em
Let $\cal D$ be the closure of the union of all forward orbits of the flow which start at the set $D$.
Under conditions above, the system has a uniquely defined, two-dimensional closed invariant set $A\subset \cal D$ $($the Lorenz attractor$)$
and a $($possibly empty$)$ one-dimensional closed invariant set $\Sigma$ $($which may intersect $A$ but is not a subset of $A)$ such that\\
$($1$)$ the separatrices $\Gamma_\pm$ and the saddle $O$ lie in $A$;\\
$($2$)$ $A$ is transitive, and saddle periodic orbits are dense in $A$;\\
$($3$)$ $A$ is the limit of a nested sequence of hyperbolic, transitive, compact invariant sets
each of which is equivalent to a suspension over a finite Markov chain with positive topological entropy;\\
$($4$)$ $A$ is structurally unstable: arbitrarily small smooth\footnote{In \cite{ABS82} this was proved for $C^1$-small perturbations,
however the argument works for $C^\infty$-small perturbations without a significant modification. A similar result is also in \cite{G76,GW}.}
perturbations of the system lead to the creation of homoclinic loops to $O$ and to the subsequent birth and/or disappearance of
saddle periodic orbits within $A$;\\
$($5$)$ when $\Sigma=\emptyset$, the set $A$ is the maximal attractor in $\cal D$;\\
$($6$)$ when the set $\Sigma$ is non-empty, it is a hyperbolic set equivalent to a suspension over a finite Markov chain
$($it may have zero entropy, e.g. be a single saddle periodic orbit$)$;\\
$($7$)$ every forward orbit in $\cal D$ tends to $A\cup \Sigma$ as $t\to+\infty$;\\
$($8$)$ when $\Sigma\neq\emptyset$, the maximal attractor in $\cal D$ is $cl(W^u(\Sigma))=A\cup W^u(\Sigma)$;\\
$($9$)$ $A$ attracts all orbits from its neighbourhood when $A\cap \Sigma =\emptyset$.\\}

Let us describe one of Shilnikov criteria for the birth of the Lorenz attractor (as defined in this theorem) from a pair of homoclinic loops.
We restrict ourselves to the symmetric case; namely, we consider a $C^r$-smooth ($r\geq 2$) system in $R^{n+1}$, $n\geq 2$, and assume
that it is invariant with respect to a certain involution $R$.

Let the system has a saddle equilibrium $O$ with $\dim W^s(O) = n-1$ and $\dim W^u(O) = 1$.
We assume that $O$ is a symmetric equilibrium, i.e. $RO = O$. With no loss of generality
we assume that the involution $R$ is linear near $O$ \cite{Bochner}.
Denote the characteristic exponents at $O$ as $\lambda_1, \dots, \lambda_n$ and $\gamma$, such that
$$\gamma>\; 0 \;> \lambda_1 > \max_{i\geq 2} {\rm Re} \lambda_i$$
(in particular, $\gamma$ and $\lambda_1$ are real). The eigenvectors $e_0$ and $e_1$ corresponding to $\gamma$ and, resp., $\lambda_1$
must be $R$-invariant. We assume that $Re_0 = -e_0$ and $Re_1 = e_1$; in particular, the two unstable separatrices are symmetric to each other,
$\Gamma_+ = R \Gamma_-$.\\

\noindent{\bf Homoclinic Butterfly.} Assume that both unstable separatrices $\Gamma_+$ and $\Gamma_-$ return to $O$ as $t \to +\infty$
and are tangent to the leading stable direction $e_1$ (it follows from the symmetry that they are tangent to each other as they enter $O$).\\

\noindent{\bf Saddle value.} Assume that the saddle value $\sigma = \gamma + \lambda_1$ is zero (hence, the saddle index $\nu = |\lambda_1|/\gamma$ equals to $1$).\\

As the homoclinic orbits $\Gamma_\pm$ are tangent to each other, we can take a common cross-section $\Pi$ to them near $O$; we can make $\Pi$ symmetric, i.e.,
$R\Pi=\Pi$. As we mentioned, the cross-section is divided by $W^s(O)$ into two parts, $\Pi_+$ and $\Pi_-$. The Poincare map $T|_{\Pi_+}$ must be symmetric
to $T|_{\Pi_-}$, so we have $u_\pm=0$, $v_+=Rv_-$, $B_+=RB_-$ and $A_+=-A_-=A$ in equation (\ref{pmps})
(the expanding direction $u$ is aligned with the vector $e_1$).\\

\noindent{\bf Separatrix value.} Assume that the separatrix value $A$ satisfies the condition $0<|A|<2$.\\

The first two conditions correspond to a codimension-2 bifurcation in the class of $R$-symmetric systems. A generic unfolding of this bifurcation
is given by a two-parameter family $X_{\mu,\varepsilon}$ of $R$-symmetric systems which depends on small parameters
$(\mu,\varepsilon)$ where $\mu$ splits the homoclinic
loops, and $\varepsilon$ is responsible for the change in the saddle value. We, thus, can
assume $\mu=u_+$ and $\nu=1-\varepsilon$ in (\ref{pmps}), and the rest of the coefficients will be continuous functions of $(\mu,\varepsilon)$. Thus, the Poincare map $T$ will have the form
\begin{equation}\label{firstret}
\begin{array}{l}
\bar u = f(u,v)=(\mu + A |u|^{1-\varepsilon} + p(|u|, Q_{s(u)} v,\mu,\varepsilon)) s(u), \\
\bar v = g(u,v)= Q_{s(u)} (v_+ + B |u|^{1-\varepsilon}) + q(|u|, Q_{s(u)} v,\mu,\varepsilon)), \\
\end{array}
\end{equation}
where $s(u)$ is the sign of $u$, and $Q_s$ denotes identity if $s=1$ and the involution $R$ (restricted to the $v$-space) if $s=-1$. The functions
$p,q$ satisfy
\begin{equation}\label{estimpq}
\left\|p,q, \frac{\partial (p,q)}{\partial v}\right\|=o(|u|^{1-\varepsilon}), \qquad \left\|\frac{\partial (p,q)}{\partial u}\right\|=o(|u|^{-\varepsilon});
\end{equation}
see \cite{book}, Section 13.8.

The following theorem is a particular case of the criteria for the birth of the Lorenz attractor which were proposed in \cite{Sh81, Sh02}.

\begin{thm} In the plane of parameters $(\mu,\varepsilon)$ there is an adjoining to zero region $V_{LA}$ such that
the system possesses a Lorenz attractor for all $(\mu,\varepsilon)\in V_{LA}$.
\label{Shil}
\end{thm}

{\bf Proof.} It is known \cite{S63} that if $\nu > 1$, then stable periodic orbits are born from homoclinic loops.
All orbits in the Lorenz attractor are unstable, so we will focus on the region $\nu < 1$, i.e.
$\varepsilon > 0$. In order to prove the theorem, we show that for all small $\varepsilon$
there exists an interval of values of $\mu$ for which the Poincare map (\ref{firstret}) satisfies conditions (\ref{abscon}) on some absorbing
domain $D$.

First of all, we check the conditions under which the invariant domain $D$ exists. According to \cite{ABS82}, the domain $D$
lies within the region bounded by the stable manifolds of saddle periodic orbits that are born from the homoclinic loops at $A\mu <0$.
Therefore, we will consider the case $\mu<0$ if $A>0$ and $\mu>0$ if $A<0$. We define
\begin{equation}\label{ddom}
D=\{|u| \leq |\mu| + \delta(\mu)\}
\end{equation}
where $\delta$ is some function of $\mu$ of order $o(\mu)$. Then, the forward invariance condition $TD\subset \mbox{int}(D)$ reads, according to
(\ref{firstret}), (\ref{estimpq}), as
$$|A| \cdot |\mu|^{1 - \varepsilon} +o(|\mu|^{1-\varepsilon}) - |\mu| < |\mu| + o(\mu),$$
or
\begin{equation} \label{ineq2}
  \displaystyle |\mu| > \left( \frac{|A|}{2}(1 + \dots)\right)^{\frac{1}{\varepsilon}}
\end{equation}
where the dots stand for the terms tending to zero as $\mu,\varepsilon\to 0$.
Note that if $|A| \geq 2$, the value at the right hand side of (\ref{ineq2}) is
bounded away from zero at small $\varepsilon > 0$, i.e., there are no small $\mu$ satisfying (\ref{ineq2}).
Thus, we consider the case $|A| < 2$, so for all
small $\varepsilon$ this inequality has solutions in any neighbourhood of the origin in the $(\mu, \varepsilon)$ plane.

Now, let us check conditions (\ref{abscon}) in $D$. By (\ref{firstret}), (\ref{estimpq}),
we have
$$f'_u=A(1-\varepsilon) |u|^{-\varepsilon}s(u)+o(|u|^{-\varepsilon}), \qquad (f,g)'_v=o(|u|^{1-\varepsilon}),\qquad
g'_u=O(|u|^{-\varepsilon}).$$
Thus,
$$
\begin{array}{ll}
\|(f'_u)^{-1}\|_{_\circ}=\frac{1}{|A|} |\mu|^{\varepsilon}(1+\dots), &   \|g'_v\|_{_\circ}=o(|\mu|^{1-\varepsilon}), \\
\|g'_u(f'_u)^{-1}\|_{_\circ}=O(1), &  \|f'_v\|_{_\circ}=O(|\mu|^{1-\varepsilon}).
\end{array}
$$
It is immediately seen from these formulas that conditions (\ref{abscon}) are fulfilled everywhere in $D$ if
\begin{equation} \label{ineq3}
\displaystyle |\mu| < \left( |A|(1 + \dots)\right)^{\frac{1}{\varepsilon}}
\end{equation}
and $\mu$ and $\varepsilon>0$ are small enough.

Conditions (\ref{ineq2}), (\ref{ineq3}) define a non-empty open region in the parameter plane near $(\mu,\varepsilon)=0$.
By construction, the intersection of this region with $\{A\mu<0, \; \varepsilon > 0\}$ corresponds to the existence of the Lorenz attractior.
$\square$

Different versions of this theorem were proved by Robinson in \cite{R92, Ro81} under an additional assumption (an open condition on the
eigenvalues of linearisation matrix at $O$). This assumption made the results of \cite{R92, Ro81} inapplicable
to the extended Lorenz model in the domain where we consider it (the case of small $\beta$ in (\ref{nfY})). We therefore
gave here a proof of the Shilnikov criterion without the extra assumptions.

The difference with the Robinson's proof is that we do a direct verfication of Afraimovich-Bykov-Shilnikov conditions for the Poincare map $T$,
while Robinson performs a factorisation of the flow by the strong-stable invariant foliation. As we mentioned, the Poincare map $T$
has a $C^1$-smooth invariant foliation $F^{ss}$. The smoothness of $F^{ss}$ is not important for
the analysis of the dynamics in the Lorenz attractor from the point of view of the topological equivalence. On the other
hand, it can be useful in the study of some statistical properties, like correlationd decay in the flow, etc. \cite{LuMe}.
The smoothness is, essentially, a consequence of two properties, that the map $T$ is hyperbolic and the expanding direction is
one-dimensional, see \cite{HPS}. Despite the hyperbolicity of the Poincare map $T$ in $D$, the suspension flow in $\cal D$
is not hyperbolic, as the set $\cal D$ contains the equilibrium state $O$. The flow is only pseudo-hyperbolic \cite{ST98,TS08},
i.e. it has a codimension-2 strong-stable invariant foliation ${\cal F}^{ss}$ such that
the flow is contracting along the leaves of the foliation and area-expanding transversely to it. The foliation $F^{ss}$ is generated by
${\cal F}^{ss}$ in the sense that the intersection of the orbit of any leaf of  ${\cal F}^{ss}$ with the cross-section $\Pi$ consists of the leaves of
$F^{ss}$. One infers from this that the foliation ${\cal F}^{ss}$ is $C^1$-smooth in ${\cal D}\backslash O$. However,
in order the foliation ${\cal F}^{ss}$ to be smooth at $O$ too,  an additional open condition on the eigenvalues of the linearisation matrix at $O$
must be fulfilled. This was the condition imposed in \cite{R92, Ro81} and the smoothness of the foliation ${\cal F}^{ss}$
was used in the Robinson proof of the birth of the Lorenz attractor in an essential way.

In what follows we apply Theorem~\ref{Shil}, in order to prove the presence of the Lorenz attractor in system
(\ref{nfY}). It suffices to find parameters $(\alpha, \lambda, \beta)$,
corresponding to the existence of a homoclinic butterfly with $\gamma=1$,  compute
the separatrix value $A$ on the homoclinig orbit, and check that $0 < A < 2$.

\section{Lorenz attractor in system (\ref{nfY}) at small $\beta$}\label{sect:main}
\subsection{Homoclinic butterfly}

At $\beta = 0$ system (\ref{nfY}) has an invariant plane $Z = 0$. If, in addition,
$\lambda = 0$, then the restriction of (\ref{nfY}) to this plane is a Hamiltonian system of the form
$\ddot X = X - X^3$. At zero energy level it posesses a pair of homoclinic solutions
$x(t) = \pm x_0(t)$, $y_0(t) = \pm \dot x_0(t)$, where
\begin{equation}\label{x0hm}
x_0(t) = \frac{\sqrt{2}}{\cosh t}.
\end{equation}
Note that by setting $\alpha = 1$ we make the saddle value $\sigma$ at the origin vanish.

Next, near the point $(\alpha, \lambda, \beta) = (1, 0, 0)$ in the space of parameters we will find the bifurcation surface which
corresponds to the existence of a homoclinic butterfly in system (\ref{nfY}). We do it by introducing a small parameter
$\mu=(\alpha-1,\lambda,\beta)$ and expanding formally the homoclinic solution in powers of $\mu$:
\begin{equation}\label{eq:loop}
X(t)=x_0(t)+x_1(t)+O(\|\mu\|^2), \qquad Z(t)=z_1(t) +O(\|\mu\|^2)
\end{equation}
where $x_1,z_1 = O(\|\mu\|)$. One obtains a linear system (with time-dependent coefficients) for the order-$\|\mu\|$ corrections
$x_1,z_1$; this system has a solution which tends to zero as $t\to \pm \infty$ when a certain linear functional of $\mu$ vanishes.
We check that this functional is not identically zero, so it vanishes on a certain plane in the $\mu$-space.
It is well known (see \cite{Kuzbook} for a proof) that this implies the existence of a surface $M_{hom}$ (tangent to this plane at $\mu = 0$)
such that system (\ref{nfY}) has a homoclinic solution of form (\ref{eq:loop}) when the parameters belong to $M_{hom}$.

We consider a more general system:
\begin{equation}
\begin{array}{l}
\ddot x + V'(x) = f(x, \dot x)z + g(x, \dot x)\mu + z^2 \phi_1(x, \dot x, z, \mu), \\
\dot z = -(\alpha(\mu) + p(x, \dot x)) z + q(x, \dot x)\mu + z^2 \phi_2(x, \dot x, z, \mu),
\end{array}
\label{nf11}
\end{equation}
where $\mu$ is an $n$--dimensional vector of small parameters, and $V$, $f$, $g$, $p$, $q$ and $\phi_{1,2}$ are certain sufficiently smooth functions.
We assume that $V(0) = V'(0) = 0$, $V''(0) < 0$, $g(0, 0) = q(0, 0) = 0$, $p(0,0)=0$, and
$\alpha(0) = \alpha_0 > 0$. By construction, the origin $O$ is an equilibrium of (\ref{nf11}) for all $\mu$.
System (\ref{nfY}) is a particular case of (\ref{nf11}), where one takes
\begin{equation}\label{vfpq}
V(x) = -\frac{x^2}{2} + \frac{x^4}{4},\;\;  f= -x, \;\; p = 0,\;\;
g = (0, -\dot x, 0), \;\; q = (0, 0, x^2).
\end{equation}
When $\mu = 0$, system (\ref{nf11}) has an invariant plane $z = 0$, and the restriction onto this plane is Hamiltonian:
\begin{equation}
\ddot x + V'(x) = 0, \\
\label{Ham}
\end{equation}
with the first integral
\begin{equation}
\dot x^2 + 2 V(x) = \mbox{const}.
\label{integr}
\end{equation}

Let there exist at least one value of $x \neq 0$ where $V(x) = 0$, and let $x = x^*$ be the closest such point to $O$
(without loss of generality we assume $x^* > 0$).
Thus, at the zero level of integral (\ref{integr}) there exists
a separatrix loop $\Gamma_0$ to $O$ lying at $z = 0$.
We take a parameterisation $x_0(t)$ of $\Gamma_0$ such that
$x_0(0) = x^*$ and, therefore, $\dot x_0(0) = 0$. It is clear that $x_0(t)$
tends to zero exponentially (with the rate $\omega = \sqrt{-V''(0)}$) as $t \to \pm \infty$.

For $\mu \neq 0$, the plane $z = 0$
is no longer invariant but the homoclinic loop $\Gamma_\mu$ may still exist if a certain
codimension-1 condition is fulfilled.
Note that if an additional $(x \to -x)$ symmetry condition is imposed on
(\ref{nf11}), then there will be two homoclinic loops to $O$.
We will seek for a homoclinic solution in form (\ref{eq:loop}).
where $(x_1(t), z_1(t)) = O(\|\mu\|)$. One obtains from (\ref{nf11}) the following system
\begin{equation}
\begin{array}{c}
\ddot x_1 + V''(x_0(t)) x_1 = \tilde f(t) z_1 + \tilde g(t) \mu,\\
\dot z_1 = -(\alpha(\mu) + \tilde p(t)) z_1 + \tilde q(t) \mu,
\end{array}
\label{linear1}
\end{equation}
where $(\tilde f, \tilde g, \tilde p, \tilde q)(t) = (f, g, p, q)(x_0(t), \dot x_0(t))$, i.e. they are explicit functions of time.

System (\ref{linear1}) is an inhomogeneous linear system; its homogeneous part is given by:
\begin{equation}
\begin{array}{c}
\ddot x_1 + V''(x_0(t)) x_1 = \tilde f(t) z_1, \\
\dot z_1 = -(\alpha(\mu) + \tilde p(t)) z_1.
\end{array}
\label{linear2}
\end{equation}
System (\ref{linear2}) has the following three linearly independent solutions:
$$(X_1(t), 0), \;\; (X_2(t), 0), \;\; (X_3(t), Z_3(t))$$
where
$$
\begin{array}{l}
X_1(t) = \dot x_0(t), \;\; X_2(t) = \displaystyle \dot x_0(t) \int\limits_0^t \frac{ds}{\dot x_0^2(s)}, \\
X_3(t) = \displaystyle  -\dot x_0(t)\int\limits_0^t \frac{dw}{\dot x_0^2(w)} \int\limits_w^{+\infty} \dot x_0(s) \tilde f(s) h(s)
e^{-\alpha s} ds, \;\;\;
Z_3(t) = e^{-\alpha t }h(t),
\end{array}
$$
and
\begin{equation}\label{eqh}
h(t) = e^{-\int\limits_0^t \tilde p(s)ds }.
\end{equation}
Note that $\dot x_0(0) = 0$ so that, generally speaking, the integrals of the form
$\int\limits_0^t \frac{F(w)}{\dot x_0^2(w)} dw$ in these formulas do not converge at $w = 0$.
We, nevertheless, keep this notation in the following exact sense:
\begin{equation}
\int\limits_0^t \frac{F(w)}{\dot x_0^2(w)} dw = \lim\limits_{\varepsilon \to +0}
\left\{ \begin{array} {l}
\displaystyle \int\limits_\varepsilon^t \frac{F(w)}{\dot x_0^2(w)} dw - \frac{F(0)}{\varepsilon \ddot x_0^2(0)}, \;\;
t > 0, \\
\displaystyle \int\limits_{-\varepsilon}^t \frac{F(w)}{\dot x_0^2(w)} dw + \frac{F(0)}{\varepsilon \ddot x_0^2(0)}, \;\;
t < 0. \\
\end{array} \right.
\label{eq:eps}
\end{equation}
While a function of the form (\ref{eq:eps}) grows to infinity as $t \to 0$, its product with
$\dot x_0(t)$ has a finite limit.

It is easy to see that the asymptotic behaviour of these solutions is represented by the following table:
$$
\begin{array}{ccc}
   & t \to -\infty & t \to +\infty \\
X_1(t) & e^{\omega t} & e^{-\omega t} \\
X_2(t) & e^{-\omega t} & e^{\omega t} \\
X_3(t), Z_3(t) & e^{-\alpha t} & e^{-\alpha t}
\end{array}
$$
where $\omega=\sqrt{-V''(0)}$.

Once we know the solution of the homogeneous system, we find the solution of the inhomogeneous system (\ref{linear1}):
\begin{equation}\label{eq:u1}
\begin{array}{c}
x_1(t) = \dot x_0(t) \displaystyle \int\limits_0^t \frac{F(w)}{\dot x_0^2(w)} dw \; \cdot \mu,  \\
\displaystyle z_1(t) = e^{-\alpha t}h(t) \int\limits_{-\infty}^t  \frac{e^{\alpha s}}{h(s)} \tilde q(s) ds \; \cdot \mu,
\end{array}
\end{equation}
where
\begin{equation}\label{fw}
F(w) =  \int\limits_{-\infty}^{w}
\dot x_0(s) \left[ \tilde f(s) h(s) e^{-\alpha s} \int\limits_{-\infty}^{s} \frac{e^{\alpha v}}{h(v)} \tilde q(v)dv +
\tilde g(s)\right] ds
\end{equation}
and $h$ is given by (\ref{eqh}).

For solution (\ref{eq:u1}) to correspond to a homoclinic orbit, it must tend to the origin in both directions of time.
It is easy to see that $\lim\limits_{t \to\pm\infty} z_1(t) = 0$ and
$\lim\limits_{t \to -\infty} x_1(t) = 0$. However, as $t \to +\infty$, the function $x_1(t)$ will converge to zero only
if
$$
F(+\infty) \cdot \mu = 0.$$
As we mentioned, when this condition defines a codimension-1 subspace in the $\mu$ space, there exists
a codimension-1 manifold $M_{hom}$ tangent to this subspace at $\mu=0$ such that system (\ref{nf11})
has a homoclinic loop at $\mu\in M_{hom}$. Thus, the surface $M_{hom}$ exists when
\begin{equation}\label{mhomc}
F(+\infty) \neq 0.
\end{equation}

Now, let us apply the obtained result to system (\ref{nfY}). By plugging (\ref{x0hm}), (\ref{vfpq}) into
(\ref{fw}) we find that
\begin{equation}\label{fw00}
\begin{array}{l}
F(w) \cdot \mu = \int\limits_{-\infty}^{w}
\dot x_0(s) \left[ \mu_3 x_0(s) e^{-s}  \int\limits_{-\infty}^{s}  x_0^2(v) e^{v} dv + \mu_2
\dot x_0(s) \right] ds = \\
= \displaystyle   -\frac{2}{3}(1 + \tanh^3 w)\lambda +
\left( \frac{55}{12} - 4 \arctan^2 e^w + \frac{15}{4 \cosh^2 w}+ \right.  \\
\left. \displaystyle +  \frac{4 (\tanh w - 2) \arctan e^w}{\cosh w}
+ 4 \tanh w + \frac{3 \tanh^2 w}{4} + \frac{4 \tanh^3 w}{3}\right) \beta.
\end{array}
\end{equation}
This gives us
$$F(+\infty) \mu= \left(\frac{32}{3} - \pi^2\right)\beta -\frac{4}{3}\lambda.$$
Condition (\ref{mhomc}) is, obviously fulfilled, so we obtain the existence of the sought homoclinic butetrfly for
parameter values belonging to a smooth $M_{hom}$ of the form
\begin{equation}
\label{eq:mu2mu3}
\beta = \frac{4}{32 - 3 \pi^2} \lambda +O(\lambda^2+(\alpha-1)^2).
\end{equation}
The corresponding homoclinic solution is given by (\ref{eq:loop}) where, according to (\ref{eq:u1}), (\ref{fw00}), we have
\begin{equation}\label{eq:x1}
\begin{array}{c}
x_1(t) = \displaystyle  -\dot x_0(t)\int\limits_0^t \frac{1}{\dot x_0^2(w)} \left[ -\frac{2}{3}(1 + \tanh^3 w)\lambda +
\left( \frac{55}{12} - 4 \arctan^2 e^w + \frac{15}{4 \cosh^2 w}+ \right. \right. \\
\left. \left. \displaystyle +  \frac{4 (\tanh w - 2) \arctan e^w}{\cosh w}
+ 4 \tanh w + \frac{3 \tanh^2 w}{4} + \frac{4 \tanh^3 w}{3}\right) \beta
\right] dw\\
z_1 = \displaystyle \mu_3 e^{-t} \int\limits_{-\infty}^t x_0^2(s)e^{s}  ds =
2e^{-t}\left( 2 \arctan e^t - \frac{1}{\cosh t} \right)\beta.
\end{array}
\end{equation}

Since the condition $\sigma = 0$ can be written as $\lambda = 2 (\alpha -1) + O((\alpha -1)^2)$, we obtain from
(\ref{eq:mu2mu3}) the following parametrisation for the bifurcation curve corresponding to the
existence of a pair of homoclinic loops to a saddle with zero saddle value in a small neighbourhood of $\mu = 0$:
\begin{equation} \label{eq:bif1}
\displaystyle
\alpha = 1 - \varepsilon + O(\varepsilon^2), \;\;
\lambda = 2\varepsilon + O(\varepsilon^2), \;\;
\beta = \frac{8}{32 - 3 \pi^2}\varepsilon + O(\varepsilon^2), \;\; \varepsilon > 0.
\end{equation}

We need to check that all conditions of the Shilnikov criterion are satisfied by the homoclinic loops that exist
when the parameters belong to this curve. First, let us check that the homoclinic loops at $\varepsilon > 0$
tend to the equilibrium along the leading direction.

Make the following change of variables in system (\ref{nfY}) in order
to align the coordinate axes with the eigendirections of the linear part of the system at $O$:
$$
X = -\alpha u + v, \;\; Y = u + \alpha v.
$$
The system takes the form
$$
\begin{array}{c}
\displaystyle \dot u = -\frac{1}{\alpha} u + \frac{\alpha}{1 + \alpha^2} uZ - \frac{1}{1 + \alpha^2} vZ - \frac{1}{1 + \alpha^2}
(v - \alpha u)^3 \\
\\
\displaystyle \dot v = \alpha v + \frac{\alpha^2}{1 + \alpha^2} uZ - \frac{\alpha}{1 + \alpha^2} vZ - \frac{\alpha}{1 + \alpha^2}
(v - \alpha u)^3 \\
\\
\dot Z = - \alpha Z + \beta (v - \alpha u)^2.
\end{array}
$$
Here the coordinates $u$ and $Z$ correspond to the stable directions, and $v$ corresponds to the unstable direction.
Note that $u$ corresponds to the strong-stable direction (i.e. the corresponding eigenvalue is farther from the imaginary axis than the eigenvalue
that corresponds to the $Z$ direction; recall that $\alpha < 1$). Thus,
the equations of the local strong stable and unstable manifolds near the saddle are
$$
\begin{array}{c}
\displaystyle
W^{ss}_{loc}: \;\; Z = \gamma_1 u^2 + O(|u|^3), \;\; v = O(u^2), \\
\displaystyle
W^{u}_{loc}: \;\; Z = \gamma_2 v^2 + O(|v|^3), \;\; u = O(v^2),
\end{array}
$$
where the coefficients $\gamma_1$ and $\gamma_2$ are found by equating the coefficients
of the power series expansion of the conditions of invariance of $W^{ss}_{loc}$ and, respectively, $W^u_{loc}$.
One obtains
$$\gamma_1 = \frac{\beta\alpha}{\alpha^2 - 2} < 0, \qquad\gamma_2 = \frac{\beta}{3\alpha} > 0.$$

This implies that the manifold
$W^{ss}_{loc}\backslash O$ lies in the region $Z < 0$, while $W^u_{loc}$ lies in the region
$Z > 0$. Finally, since $\dot Z > 0$ everywhere at $Z = 0$ except for the equilibrium point,
the global unstable manifold $W^u$ never crosses  $Z = 0$, i.e. the unstable separatrices of the saddle can never enter the region
$Z<0$, hence they cannot enter $W^{ss}_{loc}$. This proves that both homoclinic loops enter the saddle along the leading direction - the $Z$-axis;
as they come from the same side $Z > 0$, they are tangent to each other at $t=+\infty$ and  form  a homoclinic butterfly, as required by the Shilnikov
criterion.

\subsection{The separatrix value}


To verify the last condition of Theorem~\ref{Shil} we will compute the separatrix value $A$ for small values of $\beta$ and check
that it lies in the range $0 < |A| < 2$ which will finalize the proof of Theorem~\ref{thm:main}.
In order to determine the separatrix value we will use the definition from \cite{TT11}\footnote{Note that different but equivalent
definitions can be found, for example, in \cite{ABS82, R89, book}}.
Let the equation of the homoclinic loop be
$u(t) = (x(t), y(t), z(t)), \; t \in (-\infty, +\infty)$.
We consider the linearisation of (\ref{nfY}) near this solution:
\begin{equation}
\dot \xi = C(t) \xi
\label{area1}
\end{equation}
where
$$
C(t) = \left( \begin{array}{ccc} 0 & 1 & 0 \\
1 - 3 x^2(t) - z(t) & \alpha - 1 / \alpha & -x(t) \\
2 \beta x(t) & 0 & -\alpha
\end{array}\right).
$$

For any two vectors $\xi_1(t)$ and $\xi_2(t)$ satisfying (\ref{area1}) their vector product $\eta = \xi_1 \times \xi_2$
evolves by the rule
\begin{equation}
\dot \eta = D(t) \eta
\label{area2}
\end{equation}
with
$$
D(t) = \mbox{tr} C(t) \cdot \mbox{Id} - C^{\top}(t) = \left( \begin{array}{ccc} -1 / \alpha & 3 x^2(t) + z(t) - 1 & -2 \beta x(t) \\
-1 & -\alpha & 0 \\
0 & x(t) & \alpha - 1 / \alpha
\end{array}\right).
$$
As $x(t)$ and $z(t)$ tend to zero when $t \to\pm\infty$, the asymptotic behaviour of solutions of
(\ref{area2}) is determined by the limit matrix
$$
D_\infty = \left( \begin{array}{ccc} -1 / \alpha & - 1 & 0 \\
-1 & -\alpha & 0 \\
0 & 0 & \alpha - 1 / \alpha
\end{array}\right)
$$
The eigenvalues of $D_\infty$ are $0$, $-(\alpha + 1 / \alpha) < 0$ and
$\alpha - 1 / \alpha < 0$. Thus for $t \to +\infty$ each
solution of (\ref{area2}) tends to a constant times the eigenvector $\eta^* = (-\alpha, 1, 0)^{\top}$
corresponding to the zero eigenvalue.
At the backward time $t \to -\infty$ all solutions grow to infinity
except for a one-parameter family of solutions which tend to $\eta^*$ multiplied to some constant.
Thus, there is only one solution $\eta(t)$ such that $\lim\limits_{t \to -\infty}\eta(t) = \eta^*$.
For this solution, we have
$$\lim\limits_{t \to +\infty}\eta(t) = A \eta^*.
$$
The constant $A$ here is indeed the separatrix value; its absolute value is the coefficient of contraction/expansion of infinitesimal areas
near the homoclinic loop, and its sign indicates  the orientation of the loop. It is clear that
\begin{equation}
\displaystyle |A| = \sup \lim\limits_{t \to +\infty}\frac{\|\eta(t)\|}{\|\eta(-t)\|},
\label{area33}
\end{equation}
where the supremum is taken over all the solutions of (\ref{area2}).

To compute the separatrix value, perform a coordinates rotation such that
the eigendirections of the matrix $D_\infty$
become the coordinate axes. The system (\ref{area2}) takes the form:
\begin{equation}
\begin{array}{l}
\displaystyle \dot v_1 = -\frac{\alpha}{1 + \alpha^2} (3 x^2(t) + z(t)) (v_1 + v_2) +
\frac{2 \alpha \beta x(t)}{1 + \alpha^2} v_3 \\
\\
\displaystyle \dot v_2 = -\left( \alpha + \frac{1}{\alpha} \right)v_2 + \frac{\alpha}{1 + \alpha^2} (3 x^2(t) + z(t)) (v_1 + v_2)
-\frac{2 \alpha \beta x(t)}{1 + \alpha^2} v_3 \\
\\
\displaystyle \dot v_3 = \left( \alpha - \frac{1}{\alpha}\right) v_3 + x(t) v_1 + x(t) v_2.
\end{array}
\label{area3}
\end{equation}
Here $v_1$ is the coordinate in the direction of the eigenvector $\eta^*$, corresponding
to the zero eigenvalue of the limit matrix. Thus we are interested in the solution
of (\ref{area3}) satisfying the conditions $v_1(-\infty) = 1$, $v_2(-\infty) = v_3(-\infty) = 0$,
and the separatrix value is $A=v_1(+\infty)$.

Introduce a new variable $v_4 = v_1 + v_2$ to simplify the system:
\begin{equation}
\begin{array}{l}
\displaystyle \dot v_4 = -\left( \alpha + \frac{1}{\alpha} \right)v_2 \\
\\
\displaystyle \dot v_2 = -\left( \alpha + \frac{1}{\alpha} \right)v_2 + \frac{\alpha}{1 + \alpha^2} (3 x^2(t) + z(t)) v_4 -
\frac{2 \alpha \beta x(t)}{1 + \alpha^2} v_3 \\
\\
\displaystyle \dot v_3 = \left( \alpha - \frac{1}{\alpha}\right) v_3 + x(t) v_4;
\end{array}
\label{area5}
\end{equation}
the solution we are looking for satisfies $v_4(-\infty) = 1$, $v_2(-\infty) = v_3(-\infty) = 0$.
As $v_3(+\infty) = 0$ automatically, we have that $A=v_4(+\infty)$.

It is not hard to show that system (\ref{area5}) has a unique solution satisfying these conditions
for all $\beta$ and $\alpha\in(0,1]$. Just note that the problem
is equivalent to the following system of integral equations:
\begin{equation}
\begin{array}{l}
\displaystyle v_4(t) = 1 - \left( \alpha + \frac{1}{\alpha} \right) \int\limits_{-\infty}^t v_2(s) ds, \\
\displaystyle v_2(t) =  \int\limits_{-\infty}^t e^{\left( \alpha + \frac{1}{\alpha} \right) (s - t)}
\left[ \frac{\alpha}{1 + \alpha^2} (3 x^2(s) + z(s)) v_4(s) - \frac{2 \alpha \beta x(s)}{1 + \alpha^2} v_3(s)\right] ds, \\
\displaystyle v_3(t) =  \int\limits_{-\infty}^t e^{\left(\frac{1}{\alpha} -  \alpha \right) (s - t)} x(s) v_4(s) ds,
\end{array}
\label{area6}
\end{equation}
and, since $x(s)$ and $z(s)$ tend exponentially to zero as $s\to-\infty$, it immediately follows that
the integral operator on the right-hand side of this system for $t$ close to $-\infty$ is uniformly contracting
in the space of bounded continuous functions on $(-\infty,t]$ with appropriately chosen exponential weights for $v_2$ and $v_3$.
Once the solution is shown to exist up to a certain value of $t$, it is continued to all larger values of $t$ as a solution of the Cauchy
problem for system (\ref{area5}). Since the solution is obtained as a of a contracting linear operator, which depends smoothly on parameters,
the solution also depends smoothly on the parameters. Because of the uniform convergence of the integrals in (\ref{area6}) as $t\to+\infty$,
the separatrix value $A=v_4(+\infty)$ is also a smooth function of the parameters, i.e. it can be found by an asymptotic expansion.

Thus, we will seek for the solution as a power series in $\varepsilon$ using formulas (\ref{eq:bif1}).
The unknown functions are represented as
$$
v_i = v_i^0 + v_i^1 \varepsilon + O(\varepsilon^2), \;\; i = 2, 3, 4.
$$
The equation of the homoclinic loop  has the form
$$
x(t) = x_0(t) + x_1(t) + O(\varepsilon^2), \;\; y(t) = \dot x_0(t) + \dot x_1(t) + O(\varepsilon^2), \;\;
z(t) = z_1(t) + O(\varepsilon^2).
$$
For $\varepsilon = 0$, system (\ref{area5}) is rewritten as:
\begin{equation}
\displaystyle \dot v_4^0 = -2 v_2^0, \;\;
\dot v_2^0 = -2 v_2^0 + \frac{3}{2} x_0^2(t) v_4^0, \;\;
\dot v_3^0 = x_0(t) v_4^0.
\label{area7}
\end{equation}
The first two equations here do not depend on $v_3^0$ and are reduced to the equation
$$
\displaystyle \ddot v_4^0 + 2 \dot v_4^0 + \frac{3}{2} x_0^2(t) v_4^0 = 0.
$$
Its solution satisfying $v_4(-\infty) = 1$ is
\begin{equation}
\displaystyle v_4^0(t) = -\frac{\sinh t}{2\cosh^2 t} e^{-t},
\label{area9}
\end{equation}
from which we obtain
$$
\begin{array}{c}
\displaystyle
v_2^0(t) = \frac{1}{4 \cosh^2 t} (1 - 2 \tanh t), \\
\\
\displaystyle
v_3^0(t) = \frac{1}{4 \sqrt{2}} (\pi + 4 \arctan \tanh \frac{t}{2} + \frac{2}{\cosh t}(2 - \tanh t)).
\end{array}
$$

Note that $v_4^0(t)$ vanishes as $t \to +\infty$.
Therefore, the asymptotic expansion for the separatrix value has the form $A = A_1 \varepsilon + O(\varepsilon^2)$.
To complete the theorem, we need to compute $A_1$ and show that $A_1 > 0$.

The first order in $\varepsilon$ terms satisfy the following system:
\begin{equation}
\begin{array}{c}
\displaystyle \dot v_4^1 = -2 v_2^1, \;\;
\dot v_2^1 = -2 v_2^1 + \frac{3}{2} x_0^2(t) v_4^1 + f_1(t), \;\;
\dot v_3^1 = x_0(t) v_4^1 + f_2(t),\\
\\
\displaystyle f_1(t) = \frac{1}{2} (6 x_0(t) x_1(t) + z_1(t))v_4^0(t) - \frac{8}{32 - 3 \pi^2} x_0(t)v_3^0(t), \\
\\
f_2(t) = 2 v_3^0(t) + x_1(t) v_4^0.
\end{array}
\label{area8}
\end{equation}
The first two equations  do not depend on $v_3^1$, so we will solve the following equation to find $v_4^1$:
\begin{equation}
\displaystyle \ddot v_4^1 + 2 \dot v_4^1 + \frac{3}{2} x_0^2(t) v_4^1 + 2 f_1(t) = 0
\label{area10}
\end{equation}
with  the boundary conditions $v_4^1(-\infty) = 0$ and $v_4^1(+\infty) = A_1$.
Two independent solutions of the homogeneous part of (\ref{area10}) are
$$
\displaystyle y_1(t) = -\frac{\sinh t}{\cosh^2 t} e^{-t}, \;\;
y_2(t) = -\frac{\cosh^2 t + 3 t \tanh t - 3 }{2 \cosh t} e^{-t},
$$
with the Wronskian equal to
$$
\left| \begin{array}{cc}
y_1(t) & y_2(t) \\ \dot y_1(t) & \dot y_2(t)
\end{array} \right| = e^{-2t},
$$
so that the general solution of (\ref{area10}) is written as:
\begin{equation}
v_4^1(t) = C_1(t) y_1(t) + C_2(t) y_2(t),
\end{equation}
and the coefficients $C_1(t)$ and $C_2(t)$ are determined by the following formulas:
\begin{equation}
C_1(t) = 2 \int\limits_{-\infty}^t y_2(s) f_1(s) e^{2s} ds, \;\;
C_2(t) = -2 \int\limits_{-\infty}^t y_1(s) f_1(s) e^{2s} ds. \;\;
\label{area11}
\end{equation}
It remains only to calculate the limit $A_1 = \lim\limits_{t \to +\infty}(C_1(t) y_1(t) + C_2(t) y_2(t))$. We have
$$
\displaystyle
\lim\limits_{t \to +\infty}C_1(t) y_1(t) = \lim\limits_{t \to +\infty}
 \frac{\int\limits_{-\infty}^t y_2(s) f_1(s) e^{2s} ds}{\frac{1}{y_1(t)}}
 = \lim\limits_{t \to +\infty}
 \frac{y_2(t) f_1(t) e^{2t}}{-\frac{\dot y_1(t)}{y_1^2(t)}} = 0.
$$
As $\lim\limits_{t \to +\infty} y_2(t) = -1/4$, this gives us
$A_1 = 1/2 \int\limits_{-\infty}^{+\infty} y_1(t) f_1(t) e^{2t} dt$.
To compute this integral, we split it into three parts. Taking into account that
$y_1(t) x_0(t) v_4^0(t) e^{2t} = \displaystyle \frac{1}{8}(2 x_0^3(t) - x_0^5(t))$ we obtain

$$
\begin{array}{c}
\displaystyle 3 \int\limits_{-\infty}^{+\infty} y_1(t) x_0(t) x_1(t) v_4^0(t) e^{2t} dt

 = \frac{3}{8} \int\limits_{-\infty}^{+\infty} (2 x_0^3(t) - x_0^5(t)) \dot x_0(t) \int\limits_0^t \frac{F(w)}{\dot x_0^2(w)} dw
 = \\

\displaystyle \left. = \frac{3}{8} \left(\frac{x_0^4(t)}{2} - \frac{x_0^6(t)}{6} \right)
\int\limits_0^t \frac{F(w)}{\dot x_0^2(w)} dw \right|_{t = -\infty}^{t = +\infty} -
\frac{3}{48} \int\limits_{-\infty}^{+\infty} \frac{3 x_0^4(t) - x_0^6(t)}{\dot x_0^2(t)}F(t)dt =
- \frac{16}{9 (32 - 3 \pi^2)},

\end{array}
$$

and

$$
\begin{array}{c}

\displaystyle \int\limits_{-\infty}^{+\infty} y_1(t) \frac{1}{2} z_1(t) v_4^0(t) e^{2t} dt =
\frac{16}{9 (32 - 3 \pi^2)} - \frac{1}{3},\\

\displaystyle -\int\limits_{-\infty}^{+\infty} y_1(t) \frac{8}{32 - 3 \pi^2} x_0(t)v_3^0(t) e^{2t} dt =
\frac{1}{3} + \frac{2\pi^2}{32 - 3 \pi^2}.
\end{array}
$$

This result gives us the following formula for the separatrix value in a small neighborhood of point $\mu = 0$:
$$\displaystyle A = \frac{\pi^2}{32 - 3 \pi^2}\varepsilon + O(\varepsilon^2).$$

Thus, for all small $\varepsilon > 0$ we have $0 < A < 2$ which means that the orientable Lorenz
attractor is born from the homoclinic butterfly. Theorem~\ref{thm:main} is proven. $\Box$

\section{Discrete Lorenz attractors in three-dimensional diffeomorphisms}\label{sect:nf}

In this section we will apply Theorem~\ref{thm:main} to prove the birth of a discrete Lorenz attractors at certain
codimension-4 bifurcations of three-dimensional Henon-like maps.
Consider a map $(x, y, z) \mapsto (\bar x, \bar y, \bar z)$ of the form
\begin{equation}
\bar x = y,\quad \bar y = z, \quad \bar  z = B x + f(y, z, \varepsilon), \label{3dHenon}
\end{equation}
where $f$ is a smooth function, $\varepsilon$ is a set of parameters, and the Jacobian $B$ is a constant.

Fixed points of the map are given by the equations
$$
  x = y = z, \quad x(1 - B) = f(x, x).
$$
The characteristic equation at the fixed point $x = y = z = x_0$ is
$$
  \lambda^3 - A \lambda^2 - C \lambda - B = 0,
$$
where $A = f'_z(x_0, x_0)$, $C = f'_y(x_0, x_0)$. When $(A, B, C) = (-1, 1, 1)$
the fixed point has multipliers $(+1, -1, -1)$. We will study bifurcations of this point.
To do this, we shift the origin to the fixed point and introduce
small parameters $\varepsilon_1=1-B$ , $\varepsilon_2=1-C$, and $\varepsilon_3=-1-A$.
The map will take the form
\begin{equation}
\begin{array}{c}
\bar x = y, \quad \bar y = z, \quad \bar z = (1 - \varepsilon_1) x + (1 - \varepsilon_2) y - (1 + \varepsilon_3) z +
a y^2 + b yz + c z^2 + \\
+ d_1 y^3 + d_2 y^2 z + d_3 y z^2 + d_4 z^3 + \dots,
\end{array}
\label{Henon}
\end{equation}
where the dots stand for the rest of the Taylor expansion.

In \cite{GGOT13} there was shown that if the following inequality is fulfilled
\begin{equation}
  \psi = (c - a)(a - b + c) > 0,
\label{psi}
\end{equation}
then map (\ref{Henon}) near the fixed point at zero satisfies conditions proposed in \cite{SST93}, which
implies that the second iteration of the map is close, in appropriately chosen rescaled coordinates,
to the time-1 shift by the orbits of the Morioka-Shimizu system. In other words, the map is a square root of
the Poincare map for a certain small, time-periodic perturbation of the Morioka-Shimizu system. As the Morioka-Shimizu system
has the Lorenz attractor, its small time-periodic perturbation has a pseudohyperbolic attractor \cite{TS08} which corresponds
to the discrete Lorenz attractor in the root of the Poincare map. Thus, we obtain that a discrete Lorenz attractor
containing the fixed point exists in a small neighborhood of zero
for a certain region of small parameters $\varepsilon_{1,2,3}$ in
map (\ref{3dHenon}).

One can easily repeat the calculations from \cite{GGOT13} for $\psi < 0$ and obtain the same normal form in this case, but with one difference -- the scaling factor for the time will be negative, which means that the presence of the Lorenz attractor in the Shimizu-Morioka model implies
the presence of a discrete Lorenz {\em repeller} in map (\ref{3dHenon}) in the case $\psi<0$. We mentioned that
the existence of the Lorenz attractor in the Morioka-Shimizu model is definitely true but not formally proven, which means
that we are confronted with the same problem when trying to establish the existence of the discrete Lorenz attractor
or repeller in the 3D Henon map.

We bypass the problem by considering the case of an additional degeneracy; namely, we study bifurcations when
$\psi$ vanishes. The list of rescaled normal forms corresponding to a hierarchy of a certain class of degeneracies
of the $(+1,-1,-1)$-bifurcation is obtained in \cite{PST98}. It was shown there that the extended Lorenz model
does appear as a normal form in some of the degenerate cases. Therefore we can apply Theorem 1 and, thus,
obtain an analytic proof of the existence of a discrete Lorenz attractor in map (\ref{Henon}) for some
open region of parameter values; see Theorem~\ref{thmHen} below.

There are two multipliers in the formula for $\psi$, so we consider two possible cases\footnote{Note that the situation when both
$(a - b + c) = 0$ and $(a - c) = 0$ has a higher degeneracy,  i.e. we would have to deal with codimension at least five
here. We do not consider this  case in this paper.}:
\begin{equation}
\begin{array}{c}
\mbox{Case I}: \displaystyle (a - b + c) = 0, \; (a - c) \neq 0 \\
\mbox {Case II}: \displaystyle (a - b + c) \neq 0, \; (a - c) = 0.
\end{array}
\label{rho}
\end{equation}
We introduce the fourth independent small parameter $\varepsilon_4$ as follows:
\begin{equation}
\begin{array}{c}
 \mbox{Case I}: \displaystyle \varepsilon_4 = \frac{1}{4}(a - b + c) \\
 \mbox{Case II}: \displaystyle \varepsilon_4 = \frac{1}{2}(a - c).
\end{array}
\label{rho1}
\end{equation}

Consider a four-parameter family of maps of type (\ref{Henon}) which depends smoothly on $\varepsilon_{_{1, 2, 3, 4}}$.
In particular, the coefficients $a$, $b$, $c$, and $d_i$ depend on $\varepsilon_{_{1, 2, 3, 4}}$ smoothly. The following lemma
provides the normal form for this codimension-four bifurcation in both cases:

\begin{lem} \label{lemnf}
For any sufficiently small $\varepsilon_{_{1, 2, 3, 4}}$, the second iteration of map $($\ref{3dHenon}$)$ 
is close, in appropriately chosen rescaled coordinates, to the time-1 shift by the orbits of:
\begin{equation}
\begin{array}{l}
\dot X = Y, \\
\dot Y = X - \lambda Y - X Z \pm X^3 + \ldots \\
\dot Z = - \alpha Z + \beta X^2 + \ldots,
\end{array}
\label{nf1}
\end{equation}
in Case~I, and
\begin{equation}
\begin{array}{l}
\dot X = Y, \\
\dot Y = X - \lambda Y + \beta X Z + R_1 X^3 + R_2 Y Z + R_3 X Z^2 + \ldots \\
\dot Z = - \alpha Z + X^2 + Z^2 + \ldots.
\label{nf2}
\end{array}
\end{equation}
in Case~II.

Here the dots stand for vanishing at $\varepsilon = 0$ terms, the parameters
$\alpha$, $\beta$, and $\lambda$ are functions of $\varepsilon$ which can take arbitrary finite values,
and the coefficients $R_i(\varepsilon)$ are bounded as small parameters $\varepsilon$ vary.
\end{lem}

As we will show below, the sign of the $X^3$ term in formula (\ref{nf1}) coincides with the sign
of
\begin{equation} \label{G300}
  G = 4(d_1 - d_2 + d_3 - d_4) - (b - 2 c)^2.
\end{equation}
Therefore, when $G < 0$, the normal form (\ref{nf1}) for Case~I coincides with system (\ref{nfY}),  which
possesses a Lorenz attractor near $(\alpha, \beta, \lambda) = (1, 0, 0)$ according to Theorem~\ref{thm:main}.
This fact, by Lemma (\ref{lemnf}), implies the existence of a discrete Lorenz attractor in the original map (\ref{3dHenon}).
Namely, the following theorem is valid:

\begin{thm} \label{thmHen}
Assume that for some $\varepsilon = \varepsilon_0$ map $($\ref{3dHenon}$)$
possesses a fixed point $M(x_0, x_0, x_0)$ with multipliers $(1, -1, -1)$,
such that
\begin{equation}
\label{cond1}
  a - b + c = 0, \;\; G < 0,
\end{equation}
where $G$ is given by formula $($\ref{G300}$)$ and $a,b,c,d_i$ are the coefficients in the Taylor expansion $($\ref{Henon}$)$.
Then in some
small neighborhood of $\varepsilon_0$ in the parameter space there exists a domain $V_{DLA}$ such that map $($\ref{3dHenon}$)$
has a discrete Lorenz attractor when $\varepsilon \in V_{DLA}$.
\end{thm}

It remains only to prove Lemma~\ref{lemnf}.

{\bf Proof.} Let $\nu$ be the multiplier of the zero fixed point of map (\ref{Henon}) which is close to $1$:
$$
\nu = 1 - \frac{\varepsilon_1 + \varepsilon_2 + \varepsilon_3}{4}+ O(\|\varepsilon\|^2).
$$
Perform the following linear coordinate transformation:
$$
u_1 = -\nu x + y, \qquad u_2 = \nu x + (\nu - 1)y - z, \qquad u_3 = \nu(1 - \varepsilon_1)x
+ (1 - \varepsilon_1 + \nu(1 - \varepsilon_2))y + \nu^2 z.
$$
Map (\ref{Henon}) takes the following form
\begin{equation}
\begin{array}{l}
\bar u_1 = -u_1 - u_2,\\
\displaystyle \bar u_2 = -\delta_1 u_1 + (-1 + \delta_2)u_2-
\frac{1}{16}(a(2 u_1 + u_2 + u_3)^2 - b(2u_1 + u_2 + u_3) \times \\
\displaystyle
\qquad \times (2u_1 + 3u_2 - u_3) + c(2 u_1 + 3 u_2 - u_3)^2) + \frac{1}{64} (-d_1 (2 u_1 + u_2 + u_3)^3) + \\
\qquad + d_2 (2 u_1 + 3 u_2 - u_3) (2 u_1 + u_2 + u_3)^2 - d_3 (2 u_1 + 3 u_2 - u_3)^2 (2 u_1 + u_2 + u_3) + \\
\qquad + d_4 (2 u_1 + 3 u_2 - u_3)^3) + O(\|u\|^2 \|\varepsilon\| + \|u\|^4),\\

\displaystyle \bar u_3 =
(1 - \delta_3)u_3 + \frac{1}{16}(a(2u_1 + u_2 + u_3)^2 - b(2u_1 + u_2 + u_3)(2u_1 + 3u_2 - u_3) + \\
\displaystyle
\qquad + c(2u_1 + 3u_2 - u_3)^2) - \frac{1}{64} (-d_1 (2 u_1 + u_2 + u_3)^3) +
 d_2 (2 u_1 + 3 u_2 - u_3) \times \\
\qquad  \times (2 u_1 + u_2 + u_3)^2 - d_3 (2 u_1 + 3 u_2 - u_3)^2 (2 u_1 + u_2 + u_3) + \\
\qquad + d_4 (2 u_1 + 3 u_2 - u_3)^3) + O(\|u\|^2 \|\varepsilon\| + \|u\|^4),\\
\end{array}
\label{res10}
\end{equation}
where
$$
\begin{array}{l}
\displaystyle \delta_1 = \nu + \varepsilon_3 - \frac{1 - \varepsilon_1}{\nu}=
\frac{\varepsilon_1 - \varepsilon_2 + \varepsilon_3}{2} +
O(\|\varepsilon\|^2),\\ \\
\displaystyle \delta_2 = 1 - \nu - \varepsilon_3 = \frac{\varepsilon_1 + \varepsilon_2 -
3\varepsilon_3}{4} + O(\|\varepsilon\|^2),\\ \\
\displaystyle \delta_3 = 1 - \nu = \frac{\varepsilon_1 + \varepsilon_2 + \varepsilon_3}{4} +
O(\|\varepsilon\|^2).
\end{array}
$$
Note that the linear part of (\ref{res10}) is in the Jordan form at $\varepsilon = 0$.

With a close to identity polynomial change of coordinates we kill all non-resonant quadratic and cubic terms. One can see that
the system takes the following form after that:
\begin{equation}
\begin{array}{l}
\bar u_1 = -u_1 - u_2,\\
\displaystyle \bar u_2 = -\delta_1 u_1 + (-1 + \delta_2)u_2 -
\frac{1}{4}(a - c)u_1 u_3 - \frac{1}{8}(a - b - 3c)u_2 u_3 + B_{300}u_1^3 + B_{102} u_1 u_3^2 +\\
\qquad + O(\|u\|^2 \|\varepsilon\| + \|u\|^2 |u_2|+ \|u\|^4),\\

\displaystyle \bar u_3 = (1 - \delta_3)u_3 +
\frac{1}{4}(a - b + c)u_1^2 +\frac{1}{16}(a + b + c)u_3^2 + O(\|u\|^2\|\varepsilon\| +|u_1u_2|+ u_2^2 + \|u\|^3),
\end{array}
\label{res1}
\end{equation}
where
$$
\begin{array}{c}
\displaystyle B_{300} = \frac{1}{32}(-a^2 + b^2 + 5 c^2 + ab - 5 bc) + \frac{1}{8}(-d_1 + d_2 - d_3 + d_4),\\
\displaystyle B_{102} = \frac{1}{128} (2 a^2 - b^2 + ab + 2 ac - bc) +
\frac{1}{32} (-3 d_1 - d_2 + d_3 + 3 d_4).
\end{array}
$$
Next, one checks that the second iteration of the map (\ref{res1}) coincides, up to terms of order
$O(\|u\|^2 \|\varepsilon\| + \|u\|^4)$, with the time-$1$ shift by the flow of the form
\begin{equation}
\begin{array}{l}
\displaystyle \dot u_1 = \hat\rho_1 u_1 + (2 + \hat\rho_2)u_2 -
\frac{1}{4}(a - c)u_1 u_3 + \frac{1}{24} (a + 3b + 5c)u_2 u_3 + O(\|u\|^3),\\
\displaystyle \dot u_2 = \rho_1 u_1 - \rho_2 u_2 + \frac{1}{2}(a - c)u_1 u_3 -
\frac{1}{4} (b + 2c)u_2 u_3 + \tilde B_{300} u_1^3 + \tilde B_{102} u_1 u_3^2 + O(|u_2| \|u\|^2),\\
\displaystyle \dot u_3 = -\rho_3 u_3
+ \frac{1}{4} (a - b + c)u_1^2 + \frac{1}{8} (a + b + c)u_3^2  + O(|u_1u_2|+u_2^2+\|u\|^3),\\
\end{array}
\label{fl1}
\end{equation}
where $\tilde B_{300} = G / 128$,
$$\rho_3 = -\ln(1 - \delta_3) = \frac{\varepsilon_1 + \varepsilon_2 +
\varepsilon_3}{4} + O(\|\varepsilon\|^2),
$$
and
$$
\left(\begin{array}{cc} \hat\rho_1 & 1 + \hat\rho_2 \\
\rho_1 & -\rho_2 \end{array} \right)
= \ln\left(\begin{array}{cc} 1 & 1\\
\delta_1 & 1 - \delta_2\end{array}\right),
$$
so
$$
\begin{array}{l}
\displaystyle
\hat\rho_1 = -\frac{\varepsilon_1 - \varepsilon_2 + \varepsilon_3}{4} +
O(\|\varepsilon\|^2),\quad \hat\rho_2 = \frac{7\varepsilon_1 - \varepsilon_2 - 5\varepsilon_3}{24} +
O(\|\varepsilon\|^2), \\ \\
\displaystyle
\rho_1 = \frac{\varepsilon_1 - \varepsilon_2 + \varepsilon_3}{2} +
O(\|\varepsilon\|^2),\quad \rho_2 = \frac{\varepsilon_1 - \varepsilon_3}{2} + O(\|\varepsilon\|^2).
\end{array}
$$

Introduce new parameters and time via the following formulas:
$$
t = \tau / s, \; s = \sqrt{2 \rho_1}, \; \alpha = \rho_3 / s, \; \lambda = \rho_2 / s.
$$

Now we consider Cases~I and~II separately.

In Case~I with $\displaystyle \varepsilon_4 = \frac{1}{4}(a - b + c)$ we perform the following scaling of
the coordinates and parameter $\varepsilon_4$:
\begin{equation}
\displaystyle
u_1 = \nu_1 X, \; u_2 = \nu_2 Y, \; u_3 = \nu_3 Z, \; \beta = -4 \varkappa (a - c) \frac{\varepsilon_4}{s},
\label{param1}
\end{equation}
where
\begin{equation}
\displaystyle
\varkappa > 0, \; \nu_1 = (2 + \hat\rho_2)\frac{\nu_2}{s}, \; \nu_2 = \sqrt{\varkappa}s^2, \; \nu_3 = -\frac{s^2}{a - c}.
\label{param2}
\end{equation}
Then, equation (\ref{fl1}) recasts as
\begin{equation}
\begin{array}{l}
\dot X = Y+O(s), \\
\displaystyle \dot Y = X - \lambda Y - X Z + \frac {\varkappa G}{16} X^3 + O(s) \\
\dot Z = - \alpha Z + \beta X^2 + O(s).
\end{array}
\label{res2}
\end{equation}

According to \cite{PST98}, the equality $G = 0$ is the higher order  degeneracy condition, which
we do not consider here, i.e. we assume that $G \neq 0$ when $\varepsilon = 0$.
After setting $ \varkappa = 16 / G$,
system (\ref{res2}) takes the form (\ref{nf1}) up to $O(s)$-terms (recall that $\varkappa > 0$, which does not allow us
to control the sign of the $X^3$--term).

In Case~II we have $\displaystyle \varepsilon_4 = \frac{1}{2}(a - c)$ and the
scaling of the coordinates and the fourth parameter is performed as follows:
\begin{equation}
\displaystyle
u_1 = \nu_1 X, \; u_2 = \nu_2 Y, \; u_3 = \nu_3 Z, \; \beta = 2 \varkappa (a - b + c) \frac{\varepsilon_4}{s}
\label{param3}
\end{equation}
where
\begin{equation}
\displaystyle
\nu_1 = (2 + \hat\rho_2)\frac{\nu_2}{s}, \; \nu_2 = \sqrt{\varkappa}s^2, \; \nu_3 = \varkappa (a - b + c) s.
\label{param4}
\end{equation}
In this case system (\ref{fl1}) takes the form
\begin{equation}
\begin{array}{l}
\dot X = Y+O(s), \\
\displaystyle
\dot Y = X - \lambda Y + \beta X Z + \frac {\varkappa G}{16} X^3 -
\frac{1}{4}\varkappa (b + 2c)(a - b + c) Y Z + \\
\qquad + 2 \varkappa^2 (a - b + c)^2 \tilde B_{102} X Z^2 + O(s + \varepsilon_4)\\
\displaystyle
\dot Z = - \alpha Z + X^2 + \frac{1}{12}\varkappa(a - b + c)^2 Z^2 + O(s).
\label{res3}
\end{array}
\end{equation}

Finally, taking $\displaystyle \varkappa = \frac{12}{(a - b + c)^2}$ we obtain formula (\ref{nf2}). Lemma is proven. $\Box$

Now, according to Lemma~\ref{lemnf}, the flow normal form of map (\ref{3dHenon}) for the codimension-four bifurcation under consideration
coincides with system (\ref{nfY}) up to asymptotically small terms. This system has a Lorenz attractor in some domain of the parameter space
due to Theorem~\ref{thm:main}, and this implies the existence of a discrete Lorenz attractor in map (\ref{3dHenon}). Theorem~\ref{thmHen} is proven. $\Box$

We remark that system (\ref{nf2}) was not studied before, so the question of the existence of Lorenz or other attractors in this system is open. 
However, a similar system of ODEs was obtained in \cite{NOU07} as a finite-dimensional reduction of the Gray-Scott PDE model with a codimension-two
singularity. Numerical experiments with the obtained ODE model revealed a Lorenz-like chaotic behaviour \cite{NTYU07}.
Thus, the analytic study of Lorenz attractors in the normal form (\ref{nf2}) can provide results relative to the spatio-temporal chaos
in heterogeneous media.

\section*{Acknowledgements} The paper was supported by grant 14-41-00044 of Russian Scientific Foundation. The authors also acknowledge 
support by the ERC AdG grant No.339523 RGDD (Rigidity and Global Deformations in Dynamics), and the Royal Society grant IE141468.

\end{document}